\documentclass{article}

\usepackage[english]{babel}

\usepackage[letterpaper,top=2cm,bottom=2cm,left=3cm,right=3cm,marginparwidth=1.75cm]{geometry}

\usepackage{amsmath,amsfonts,amssymb, stmaryrd,mathtools,amsthm}
\usepackage{graphicx}
\usepackage[colorlinks=true, allcolors=blue]{hyperref}
\allowdisplaybreaks[1]

\newtheorem{theorem}{Theorem}[section]

\newtheorem{lemma}[theorem]{Lemma}

\newtheorem{remark}[theorem]{Remark}

\newcommand{\q}{_{B^{2-2/p}_{q,p}(\mathbb{T})}}
\newcommand{\p}{_{L^q(\mathbb{T})}}
\newcommand{\f}{_{\infty}}
\newcommand{\T}{\mathbb{T}}
\newcommand{\R}{\mathbb{R}}
\newcommand{\fl}{\int_{\mathbb{T}}f\phi_ldr}
\newcommand{\1}{\partial_{r}}
\newcommand{\2}{\partial_{rr}}

\title{Stochastic Curve Shortening Flow with Scale-Dependent Noise }
\author{Qi Yan\footnote{Academy of Mathematics and Systems Science, Chinese Academy of Sciences, Zhuangguancun East Road 55, Beijing 100190, China,  and Universität Leipzig, Fakultät für Mathematik und Informatik, Augustusplatz 10, 04109 Leipzig, Germany}}

\begin{document}
\maketitle

\begin{abstract}
In this paper, we study the motion by mean curvature of curves in the plane perturbed by scale-dependent noise. We first introduce a so-called scale-dependent noise from the physics background to the curve shortening flow. To be more precise, the scale-dependent noise defined on a curve is a noise whose intensity is proportional to the length of the curve. To get the well-posedness of stochastic curve shortening flow driven by scale-dependent noise, we equivalently formulate 
the stochastic curve shortening flow as a one-phase stochastic Stefan problem of its curvature parameterized by the arclength parameter and its length. After rewriting the one-phase stochastic Stefan problem as a quasilinear evolution equation, we apply the theory for quaslinear stochastic evolution equations developed by Agresti and Veraar in 2022 to get  maximal unique local strong solution for the stochastic curve shortening flow up to a maximal stopping time which is characterized by a blow-up criterion.
\\\\
{\bf Ketwords:} Stochastic curve shortening flow,  Stochastic one-phase Stefan problem, scale-dependent noise
\end{abstract}

\section{Introduction}
Mean curvature flow  as following of embedded hypersurfaces in Euclidean spaces is a very important geometric flow and has been intensively studied in the past decades, see \cite{Hui84,Hui90,CM12}. 
\[
    \frac{\partial F}{\partial t}=-H\mathbf{n}
\]
where $F:M^n\times [0,T)\to\mathbb{R}^{n+1}$ is an embedding map, and $H$ is the mean curvature, $\mathbf{n}$ is the unit normal vector.
Due to its gradient flow structure (mean curvature flow is the $L^2$-gradient flow of the area functional of hypersurfaces)  mean curvature flow describes how interfaces evolve to reduce their surface energy, leading to a gradual smoothing or flattening of the surface.  In physics, it arises as the sharp interface limit of Allen-Cahn equation for the phase field of a binary alloy, describing the motion of an interface between two phases. Stochastic mean curvature flow was first proposed in \cite{KO82} as a refined model incorporating the influence of thermal noise. As a result one may think of the following stochastic evolution
\begin{equation}
    dF_t=-\left(Hdt+{W}(M_t,x,\circ dt)\right)\mathbf{n},\quad x\in M_t
\end{equation}
where $M_t=F_t(M)$, $W$ is a specific random field and $W(M_t,x,\circ dt)$ is its Stratonovich differential. As an example, we can choose $W(M_t,x,t)=\phi(M_t,x)b_t$ for some smooth function $\phi$ and a standard one-dimensional Brownian motion $b_t$, inducing the dynamics
\begin{equation}
    dF_t=-\left(Hdt+\phi(M_t,x)\circ db_t\right)\mathbf{n}
\end{equation}

Mean curvature flow on the plane is also called the curve shortening flow, 
\begin{equation}
    \partial_t\gamma=-k\mathbf{n}
\end{equation}
where $\gamma:\mathbb{S}^1\times[0,T)\to\mathbb{R}^2$, $k$ is the curvature of the curve $\gamma_t$ and $\mathbf{n}$ is the normal vector. A well-known result of the curve shortening flow is that for any simple closed curve $\gamma_0$ on the plane, it will eventually shrink to a round point in the plane in finite time, during it shrinks, it becomes convex at some time and remain so, see \cite{GH86, Gra87}.

In this paper we will consider the following stochastic curve shortening flow
\begin{equation}\label{Stochastic curve shortening flow}
    d \gamma_t =-\left(kdt+\sigma L(t)\circ dW_t\right)\mathbf{n}
\end{equation}
where $L(t)$ is the length of the curve $\gamma_t$, $\sigma>0$ is a constant and $W_t$ is a standard one-dimensional Brownian motion.  We would like to call the noise $\sigma L(t)\circ dW_t$ a scale-dependent noise.

We next explain the motivation for introducing such noise.  

In many physical systems, noise is generated by the superposition of multiple independent random sources in the system. The larger the size of the system, the more independent noise sources it contains, so the total noise intensity is proportional to the size. For example, Thermal noise in the wire in electronics, also called Johnson noise, the thermal noise voltage  of the wire is directly proportional to the length $L$ of the wire. For a curve, the length $L$ is its natural dimensional parameter, so it is reasonable to assume the noise intensity depends on its length. 

If we place a single particle in a noise environment, like a grain of pollen in water as Robert Brown did, Einstein in 1905 in \cite{Einstein1905} proved that the dynamic of the particle is actually a Brownian motion $\lambda W_t$, where $\lambda$ is a constant depending on the mass, velocity, radius of the particle, etc., and $W_t$ is the standard Brownian motion on the plane. Now, we consider a curve shortening flow on the plane, every point on the curve receives infinite collision from particles in the environment. Since collision along the tangent direction does not change the shape of the curve, we only consider collision along the normal direction. According to Einstein's theory, the displacement along the normal direction caused by collision in a short time $dt$ of a point is $\sigma dW_t$, where $\sigma$ is another constant and $W_t$ is a 1-dimensional Brownian motion. Since the displacement of a point on a curve will affect its neighborhoods on the curve, the displacement along the normal direction of a point on a curve should be the displacement accumulation along the curve, hence, the displacement along the normal direction caused by the collision at every point on a curve with length $L$ should be $\sigma L\circ dW_t$. Thus, the displacement along the normal direction at each point of the curve shortening flow in a short time $dt$ should be $kdt+\sigma L(t)\circ dW_t$. Thus, we can say that (\ref{Stochastic curve shortening flow}) models the curve shortening flow in a noise environment.

Stochastic curve shortening flow perturbed by one dimensional Brownian motion and infinite dimensional Brownian motion have been studied by Es-Sarhir and von Renesse in \cite{EvR12}, where they studied the stochastic curve shortening flow of curves which can be expressed as the graph of some function defined on one-dimensional torus $\mathbb{T}$, to be more precise,  the SPDE they studied are as following:
\begin{equation}
    du(x)=\frac{\partial^2_xu}{1+(\partial_xu)^2}(x)dt+\sqrt{1+(\partial_xu)^2 }\circ dW.
\end{equation}
the stochastic evolution equation of the some functions whose graphs as curves on $\mathbb{R}^2$ satisfy the stochastic curve shortening flow. Following this paper, stochastic mean curvature flows in $\mathbb{R}^3$ as well as $\mathbb{R}^n$ were studied in \cite{HRvR17} and \cite{DHR21} respectively. However, the way of studying the evolution equations of the graph functions of general curve flow of curves in the plane cannot be applied to general case, for example, closed or non-periodic planar curves cannot be expressed as graphs of some functions defined on $\mathbb{T}$. 

Thus, we would like to formulate the stochastic curve shortening flow \eqref{Stochastic curve shortening flow} in a different way. The evolution equations we study in this paper are as following:
\begin{equation}\label{evolution equation system}
    \left\{
    \begin{aligned}
        dk(t)&=\left(\partial_{ss}k+k^3\right)dt+\sigma k^2L(t)\circ dW_t,\quad s\in [0,L(t)],\\
        dL(t)&=-\int^{L(t)}_0k^2dsdt-\sigma L(t)\int^{L(t)}_0kds\circ dW_t,\\
        k(s,0)&=k_0(s), \,L(0)=L_0,\,s\in [0,L_0].
    \end{aligned}
    \right.
\end{equation}
where $k(t)$ and $L(t)$ are the curvature and length of the curve $\gamma_t$, $s$ is the arclength parameter of $\gamma_t$, and $\partial_{ss}$ denotes the Laplace-Beltrami operator of the curve $\gamma_t$.

In fact, \eqref{evolution equation system} is a stochastic one-phase Stefan problem (we can also call it a stochastic free boundary problem or a stochastic moving boundary problem). As long as we get a solution $(k(s, t), L(t))$, we can construct the curve $\gamma_t$ in the following way:
\begin{equation*}
    \begin{aligned}
    \gamma_t(s)=\left(x_0+\int^s_0\cos\left(\theta_0+\int^r_0k(u)du\right)dr,\,y_0+\int^s_0\sin\left(\theta_0+\int^r_0k(u)du\right)dr\right),\quad s\in[0,L(t)]
    \end{aligned}
\end{equation*}
where $(x_0,y_0)$ and $\theta_0$ are decided by the starting point $\gamma_0(0)$ and tangent vector $\gamma^{\prime}_0(0)$.

Compared to the way of studying the evolution equations of the graph functions of a general curve
flow, our approach of studying the evolution equations of the curvature and length of a curve flow is
more universal, since any initial curve, no matter closed or not, with curvature can be studied by this
approach. In addition, all quantities in \eqref{evolution equation system} are intrinsic, that is to say they do not depend on
the ambient space where the curves lie in.

\section{Preliminaries}
In this section, an overview over the basic tools and notions used in this paper is provided. For more details, we give references to the literature.

Throughout this paper, we fix a probability space $(\Omega,\mathcal{F},\mathbb{P})$ with filtration $\{\mathcal{F}_t\}_{t\geq 0}$, a separable Hilbert space $H$, which satisfies the usual conditions. Moreover, for two given normed spaces $X$ and $Y$, the set of all linear operators from $X$ to $Y$ is denoted by $B(X,Y)$.
\subsection{Stochastic integration}
\subsubsection{The space $\gamma(\mathcal{H},X)$}
Let $\mathcal{H}$ be a Hilbert space(typically, we take $\mathcal{H}=H$ or $\mathcal{H}=L^2(0,T;H)$). The Banach space $\gamma(\mathcal{H},X)$ of all $\gamma$-radonifying operators from $\mathcal{H}$ to $X$ defined as the closure of the space of finite rank operators from $\mathcal{H}$ to $X$ with respect to the following norm
\[
\|T\|^2_{\gamma(\mathcal{H},X)}:=\mathbb{E}\left\|\sum^{\infty}_{n=1}\gamma_nTh_n\right\|^2_X,
\]
where $(h_n)_{n\in\mathbb{N}}$ is an orthonormal basis of $\mathcal{H}$, and $(\gamma_n)_{n\in\mathbb{N}}$ is any sequence of independent standard Gaussian random variables. Note that the norm is independent of the choice of the orthonormal basis. To know more about the theory of $\gamma$-radonifying operators, we refer to \cite{DJT} and \cite{vanN}. 

In the special case, that $X=L^p(O,\mu)$ with $p\in[1,\infty)$ and $(O,\mu)$ $\sigma$-finite, one has a canonical isomorphism:
\[
L^p(O,\mu;\mathcal{H})\simeq\gamma(\mathcal{H},L^p(O,\mu)),
\]
which is obtained by the following mapping $L^q(O;\mathcal{H})\ni f\mapsto T_f\in \gamma(\mathcal{H};X)$, where $T_f$ is defined by
\[
T_f(h)(x):=\langle f(x),h\rangle_{\mathcal{H}}\quad \forall h\in \mathcal{H} \text{   and  } x\in O.
\]
The equivalence $\|T_f\|_{\gamma(\mathcal{H};X)}\simeq\|f\|_{L^q(O;\mathcal{H})}$ can be shown easily by the Kahane-Khintche inequality 
\begin{equation}
    \left(\mathbb{E}\left\|\sum_{l\in\mathbb{N}}\gamma_lf_l\right\|^q_X\right)^{1/q}\simeq_q\mathbb{E}\left\|\sum_{l\in\mathbb{N}}\gamma_lf_l\right\|_X,
\end{equation}
for $q\in[1,\infty)$.

\subsubsection{The stochastic integral}
For a stochastic process $G:\Omega\times\mathbb{R}_+\times H\to X$ of the form
\[
G=\chi_{(s,t]\times F}h\otimes x,
\]
with $F\in\mathcal{F}_s, h\in H$ and $x\in X$, we can define the stochastic integral via 
\[
I(G):=\int^T_0GdW:=\chi_FW(\chi_{(s,t]}h)x,
\]
 and then extend it to $\mathcal{F}$-adapted step processes, which are finite linear combinations of such processes. Van Nerven, Veraar, and Weis proved in \cite{vNVW07} the following two-sided estimate for this stochastic integral.
 \begin{theorem}
     Let $X$ be a UMD Banach space, and $G$ be an $\mathcal{F}$-adapted step process in $\gamma(H;X)$. Then for all $p\in (0,\infty)$, one has the two-sided estimate
     \[
     \mathbb{E}\|I(G)\|^p_X\simeq_p\mathbb{E}\|G\|^p_{\gamma(L^2(0,T;H);X)},
     \]
     with implicit constants depending only on $p$ and the UMD constant of $X$.
     In particular, the stochastic integral can be continued to a linear and bounded operator
     \[
     I:L^p(\Omega;\gamma(L^2(0,T;H);X))\longrightarrow L^p(\Omega;X).
     \]
 \end{theorem}
 All Hilbert spaces and Banach spaces $L^q(D,\mu)$ with $q\in (1,\infty)$ are UMD spaces. Furthermore, closed subsets, quotients and duals of UMD spaces are UMD. For a UMD Banach space $X$ with type 2(more details of UMD spaces and type of them can be found in \cite{BURKHOLDER2001233} and \cite{PG}), for $q>2$, we have the following continuous embedding
 \[
 L^p(0,T;\gamma(H;X))\hookrightarrow L^2(0,T;\gamma(H;X))\hookrightarrow\gamma(L^2(0,T;H);X).
 \]
  Thus, the stochastic integral $I(G)$ can be defined for $G\in L^p(\Omega\times[0,T];\gamma(H;X))$.

  \subsection{$R$-boundedness and $H^{\infty}$-calculus}
  Let $X$ and $Y$ be two Banach spaces and $(r_n)_{n\geq 1}$ be a sequence of Randemacher random variables, i.e. $\mathbb{P}(r_n=1)=\mathbb{P}(r_n=-1)=1/2.$ A family $\mathcal{T}\subset B(X,Y)$ is called $R$-bounded, if there is a constant $C>0$, such that
  \[
  \mathbb{E}\left\|\sum^N_{n=1}r_nT_nx_n\right\|^2_Y\leq \mathbb{E}\left\|\sum^N_{n=1}r_nx_n\right\|^2_X,
  \]
  for any finite sequence $(T_n)^N_{n=1}\subset \mathcal{T}$ and $(x_n)^N_{n=1}\subset X$. The least admissible constant $C$ is called the $R$-bound of $\mathcal{T}$, denoted by $R(\mathcal{T})$. Note that every $R$-bounded family is uniformly bounded family in $B(X,Y)$.If $X$ and $Y$ are Hilbert spaces, the $R$-boundedness is equivalent with uniform boundedness and $R(\mathcal{T})=\sup_{T\in \mathcal{T}}\|T\|$. For more details of $R$-boundedness and its application, we refer to \cite{DHP,CPS,KW}

  An operator $A$ with domain $D(A)$ is called sectorial on a Banach space $X$ of angle $\theta\in (0,\pi/2)$, if it is closed, densely defined, injective and has a dense range, moreover, its spectrum is contained in the sector $\Sigma_{\theta}=\{z\in \mathbb{C}:|\arg(z)|<\theta\}$ and the set \[
  \{\lambda R(\lambda,A):\lambda\notin \Sigma_{\phi}\},
  \]
  is bounded in $B(X)$ for all $\phi\in (\theta,\pi)$ and the bound only depends on $\phi$. In this case, $-A$ generates a holomorphic semigroup on $X$.

  For any holomorphic function $f$ on $\Sigma_{\phi}, \phi\in (\theta,\pi)$, satisfying the estimate
  \[
  |f(z)|\leq C\frac{|z|^{\delta}}{1+|z|^{2\delta}},
  \] for some $\delta>0$, (the space of those function is denoted by $H^{\infty}_0(\Sigma_{\phi})$), the integral
  \[
  f(A)=\frac{1}{2\pi i}\int_{\partial\Sigma_{\phi}}f(z)R(z,A)dz,
  \]
  converges absolutely and is independent of $\phi$. We say that $A$ has a bounded $H^{\infty}(\Sigma_{\theta})$-calculus, if there exists a constant $C>0$ such that 
  \[
  \|f(A)\|_{B(X)}\leq C\|f\|_{\infty},\quad\quad \forall f\in H^{\infty}_0(\Sigma_{\theta}).
  \]
  The least constant $C$ is called the bound of $H^{\infty}$-calculus.

\section{Equivalence of a General curve flow and a one-phase Stefan problem}
In this section, we prove that a general curve flow is equivalent with a one-phase moving boundary problem of its curvature and length.

For a general curve flow 
\begin{equation}\label{general flow}
    \frac{\partial F}{\partial t}=-VN,\quad F(0)=F_0
\end{equation}
where $F:S^1\times [0,T)\to \mathbb{R}^2$ represent a one parameter family of closed curves with counterclockwise parameterization, $N$ is the outward-pointing unit normal, $V$ is the shrinking speed in the normal direction. Next, we will derive the evolution equations for its curvature $k$ and the length of the curves $L$.

First, we reparameterize the curve $F(u,t)$ by its arclength parameter $s$. The arclength parameter $s$ is defined by \[
s=s(u):=\int^u_0\left|\frac{\partial F}{\partial r}(r,t)\right|dr
\]
we denote $v=|\partial F/\partial u|$, then we have $ds=vdu$ and \[ \frac{\partial}{\partial s}=\frac{1}{v}\frac{\partial}{\partial u}. \]
Let $T$ be the unit tangent vector to the curve $F(u)$. The Frenet-Serret equations are 
\[
\frac{\partial T}{\partial s}=-kN,\quad \frac{\partial N}{\partial s}=kT,
\]
 in terms of $u$, it becomes 
 \[
 \frac{\partial T}{\partial u}=-vkN,\quad \frac{\partial N}{\partial u}=vkT
 \]
 Note that
 \begin{align*}
     \frac{\partial }{\partial t}\left(v^2\right)&=\frac{\partial}{\partial t}\left\langle\frac{\partial F}{\partial u},\frac{\partial F}{\partial u}\right\rangle=2\left\langle\frac{\partial F}{\partial u},\frac{\partial^2 F}{\partial t\partial u}\right\rangle=2\left\langle\frac{\partial F}{\partial u},\frac{\partial^2 F}{\partial u\partial t}\right\rangle=2\left\langle vT, \frac{\partial }{\partial u}(-VN)\right\rangle\\
     &=2\left\langle vT, \,-\frac{\partial V}{\partial u}N-VvkT\right\rangle=-2Vv^2kN
 \end{align*}
 Note here $\partial/\partial u$ and $\partial/\partial t$ commute since $u$ and $t$ are independent coordinates. Hence, we have 
 \[
 \frac{\partial v}{\partial t}=-Vkv
 \]
 Then, we can easily get \[
 \frac{\partial L}{\partial t}=\frac{\partial}{\partial t}\int^{2\pi}_0vdu=\int^{2\pi}_0\frac{\partial v}{\partial t}du
=-\int^{2\pi}_0Vkvdu=-\int^L_0Vkds \]

Before we proceed to get the evolution equation of the curvature $k$, we need the following relation for the operators $\partial /\partial s$ and $\partial/\partial t$,
\[ \frac{\partial}{\partial t}\frac{\partial}{\partial s}=\frac{\partial}{\partial t}\frac{1}{v}\frac{\partial}{\partial u}=Vk\frac{1}{v}\frac{\partial}{\partial u}+\frac{1}{v}\frac{\partial}{\partial u}\frac{\partial}{\partial t}=\frac{\partial}{\partial s}\frac{\partial}{\partial t}+Vk\frac{\partial}{\partial s}. \]
Next, we need the time derivatives of $T$ and $N$.
\begin{align*}
    \frac{\partial T}{\partial t}&=\frac{\partial^2F}{\partial t\partial s}=\frac{\partial^2F}{\partial s\partial t}+Vk\frac{\partial F}{\partial s}=\frac{\partial}{\partial s}\left(-VN\right)+VkT=-\frac{\partial f}{\partial s}N-VkT+VkT=-\frac{\partial V}{\partial s}N
\end{align*}
 since we have 
 \[
 0=\frac{\partial}{\partial t}\left\langle T,N\right\rangle=\left\langle\frac{\partial V}{\partial t}N, N\right\rangle+\left\langle T,\frac{\partial N}{\partial t}\right\rangle
 \]
 we get
 \[
 \frac{\partial N}{\partial t}=\frac{\partial V}{\partial s}T
 \]

Let $\theta$ be the angle between the tangent vector $T$ and the $x$-axis, we can write $T=(\cos\theta,\sin\theta),N=(\sin\theta,-\cos\theta)$, then
\[
\frac{\partial T}{\partial t}=-\frac{\partial V}{\partial s}N=-\frac{\partial V}{\partial s}(\sin\theta,-\cos\theta),\quad \frac{\partial T}{\partial s}=-kN=-k(\sin\theta,-\cos\theta)
\]
on the other hand,
\[
\frac{\partial T}{\partial t}=\frac{\partial}{\partial t}(\cos\theta,\sin\theta)=\frac{\partial\theta}{\partial t}(-\sin\theta,\cos\theta),\quad \frac{\partial T}{\partial s}=\frac{\partial }{\partial s}(\cos\theta,\sin\theta)=\frac{\partial\theta}{\partial s}(-\sin\theta,\cos\theta)
\]
 Thus, we have \[
 \frac{\partial\theta}{\partial t}=\frac{\partial V}{\partial s},\quad \frac{\partial\theta}{\partial s}=k
 \]
 
 The evolution equation for the curvature $k$ will be \[
 \frac{\partial k}{\partial t}=\frac{\partial^2\theta}{\partial t\partial s}=\frac{\partial^2 \theta}{\partial s\partial t}+Vk\frac{\partial\theta}{\partial s}=\frac{\partial^2 V}{\partial s^2}+k^2 V,\quad s\in[0,L(t)].
 \]

Since every planar curve is uniquely determined by its curvature $k(s,t),s\in [0,L(t)]$ up to translation and rotation in $\mathbb{R}^2$. As long as we know the curvature function $k(s)$ of a curve $\gamma(s) $ parameterized by the arc-length parameter $s\in [0,L]$, we can construct the curve as follows:

The angle between the tangent vector and $x$-axis is
\[
\theta(s)=\theta_0+\int^s_0k(r)dr
\]
then the tangent vector $T=(\cos\theta,\sin\theta)$ is 
\[
T(s)=\left(\cos\theta(s),\sin\theta(s)\right)=\left(\cos\left(\theta_0+\int^s_0k(r)dr\right),\,\sin\left(\theta_0+\int^s_0k(r)dr\right)\right)
\]
Thus, the curve $\gamma(s)$ is
\begin{equation}
    \begin{aligned}
    \gamma(s)&=\gamma_0+\int^s_0T(r)dr\\
        &=\left(x_0+\int^s_0\cos\left(\theta_0+\int^r_0k(u)du\right)dr,\,y_0+\int^s_0\sin\left(\theta_0+\int^r_0k(u)du\right)dr\right)
    \end{aligned}
\end{equation}

Hence, the general curve flow  (\ref{general flow}) is equivalent with the following one-phase Stefan problem of the curvature $k(s,t)$ and length $L(t)$
\begin{equation}\label{general evolution equation}
    \left\{\begin{aligned}
    &\partial_t k=\partial_{ss}V+k^2 V,\quad s\in [0,L(t)]\\
    &\partial_t L(t)=-\int^{L(t)}_0kVds\\
    &k(s,0)=k_0(s), \quad L(0)=L_0, \quad s\in [0,L_0]
\end{aligned}
\right.
\end{equation}
Note that here the curvature $k(s,t)$ of closed curves should be understood as an $L(t)$-periodic function.

\section{Existence of Local Solution of Stochastic curve shortening flow}
The shrinking speed in the normal direction for our stochastic curve shortening flow (\ref{Stochastic curve shortening flow}) is
\[
V=k+\sigma L(t)\circ dW_t
\]
By (\ref{general evolution equation}), the evolution equations of the curvature $k$ and length $L$ of the curve are 
\begin{equation}
    \left\{
    \begin{aligned}
        dk(t)&=\left(\partial_{ss}k+k^3\right)dt+\sigma k^2L(t)\circ dW_t,\quad s\in [0,L(t)]\\
        dL(t)&=-\int^{L(t)}_0k^2dsdt-\sigma L(t)\int^{L(t)}_0kds\circ dW_t\\
        k(s,0)&=k_0(s), \,L(0)=L_0,\,s\in [0,L_0]
    \end{aligned}
    \right.
\end{equation}
Note that, for any simple close planar curve $\gamma$, we have 
\[
\int_{\gamma}kds=2\pi
\]
Thus, the evolution equations for closed curves are 
\begin{equation}
    \left\{
    \begin{aligned}
        dk(t)&=\left(\partial_{ss}k+k^3\right)dt+\sigma k^2L(t)\circ dW_t,\quad s\in [0,L(t)]\\
        dL(t)&=-\int^{L(t)}_0k^2dsdt-2\sigma\pi L(t)\circ dW_t\\
        k(s,0)&=k_0(s), \,L(0)=L_0,\,s\in [0,L_0]
    \end{aligned}
    \right.
\end{equation}
We make the following transform\[
s=rL(t), \quad f(r,t)=k(rL(t),t),\quad r\in \mathbb{T}
\]
then the evolution equation of $L$ becomes
\[
dL(t)=-L\int_{\mathbb{T}}f^2dr-2\sigma\pi L\circ dW_t
\]
since the Stratonovich differential satisfies the chain rule, we have 
\begin{align*}
    df(t)=\left(\frac{1}{L^2}\2f+f^3-r\1f\int_{\T} f^2dr\right)dt+\sigma\left(f^2L-2\pi r\1f\right)\circ dW_t
\end{align*}
Next, we transform the Stratonovich differential into It\^{o} differential.
\[
L\circ dW_t=LdW_t+\frac{1}{2}\langle dL,dW\rangle=LdW_t-\sigma \pi Ldt
\]
\begin{align*}
    \1f\circ dW_t&=\1fdW_t+\frac{1}{2}\langle d\1f,dW\rangle=\1f dW_t+\frac{1}{2}\langle\partial_r df,dW\rangle\\
    &=\1f dW_t+\left(\sigma f\1fL -\sigma \pi \1f-\sigma \pi r \2f\right)dt
\end{align*}
\begin{align*}
    f^2L\circ dW_t&=f^2LdW_t+\frac{1}{2}\langle d(f^2L),dW\rangle=f^2LdW_t+fL\langle df,dW\rangle+\frac{f^2}{2}\langle dL,dW\rangle\\
    &=f^2LdW_t+\left(\sigma f^3L^2 -2\sigma\pi rf\1fL-\sigma\pi f^2 L\right)dt
\end{align*}
Thus the evolution equations of $(f,L)$ are
\begin{equation}\label{equation for (f,L)}
    \left\{
    \begin{aligned}
        df=&\left[\left(2\sigma^2\pi^2r^2+\frac{1}{L^2}\right)\2f-4\sigma^2\pi rLf\1f+2\sigma^2\pi^2 r \1f -r\1f\int_{\T} f^2dr+f^3+\sigma^2f^3L^2-\sigma^2\pi f^2L\right]dt\\
        &\quad +\sigma\left(f^2L-2\pi r\1f\right)dW_t\\
        dL=&L\left(2\sigma^2\pi^2-\int_{\T} f^2dr\right)dt -2\sigma\pi LdW_t\\
        f(r,0)&=k_0(rL_0),\quad L(0)=L_0, \quad r\in\mathbb{T}
    \end{aligned}
    \right.
\end{equation}
\begin{remark}
    In fact, we can directly get a precise expression for $L$ in term of $f$
    \begin{equation}
        L(t)=L_0\exp\left(-\int^t_0\int_{\T} f^2(r,\tau)drd\tau-2\sigma\pi W_t\right)
    \end{equation}
    which is similar with the geometric Brownian motion.
\end{remark}

We define the linear operator $A(f,L)$ by 
\[
A(f,L)=-\left[\left(2\sigma^2\pi^2r^2+\frac{1}{L^2}\right)\2\right]I_2,
\]
where $I_2=\mathrm{diag}(1,1)$, the 2-dim identity matrix,
define the nonlinear term$F(t,f,L)=(F_1(t,f,L), F_2(t,f,L))^{\top}$, where
\begin{align*}
    F_1(t,f,L)&:=-r\1f\int_{\T}f^2dr+f^3+\sigma^2\left(-4\pi rLf\1f+2\pi^2 r \1f +f^3L^2-\pi f^2L\right), \\
    F_2(t,f,L)&:=2\sigma^2\pi^2 L-L\int_{\T}f^2dr.
\end{align*}
define $B(t,f,L)=B_1(t,f,L)+B_2(t,f,L)$ where
\begin{equation*}
    B_1(f,L):=\left(\begin{array}{cc}
         \sigma f^2L  \\
          0
    \end{array}\right),
    \quad\quad 
    B_2(f,L):=\left(\begin{array}{cc}
         -2\sigma\pi r\1f  \\
          -2\sigma\pi L
    \end{array}\right)
\end{equation*}
With the above notations, we can rewrite the evolution equation (\ref{equation for (f,L)}) of $(f,L)$ in the following quasilinear form
\begin{equation}\label{quasilinear equation}
    \begin{aligned}
        d(f,L)^{\top}&=\left[-A(f,L)(f,L)^{\top}+F(t,f,L)\right]dt+B(f,L)dW_t,\\
        (f(0),L(0))&=(f_0,L_0).
    \end{aligned}
\end{equation}

We choose the spaces\footnote{In fact, here we should choose $X=L^{q}(\mathbb{T})\times L^q(\mathbb{T})$, but in our case, the function $L(t)$ is only a real function of time $t$, which means $L(t) \in \mathbb{R},\|L\|\q=|L|$, so we just choose $X=L^q(\mathbb{T})\times\mathbb{R}$ for simplicity, so does $X_1$.}$X=L^q(\mathbb{T})\times \mathbb{R}, X_1=W^{2,q}(\mathbb{T})\times \mathbb{R}$ and their norms are defined by 
\begin{equation*}
    \|(f,L)\|_X=\|f\|_{\q}+|L|\quad \|(g,M)\|_{X_1}=\|g\|_{W^{2,q}(\mathbb{T})}+|M|
\end{equation*}
for all $(f,L)\in X$ and $(g,M)\in X_1$. 
Since $L^q(\mathbb{T}), W^{2,q}(\mathbb{T})$ and $\mathbb{R}$ are all UMD Banach spaces with type 2, it is easy to check that the product spaces $X, X_1$ with the above norms are also UMD Banach spaces with type 2. What is more, the real interpolation space \begin{align*}
    X_p&=(X,X_1)_{1-1/p,p}=(L^q(\mathbb{T})\times \mathbb{R},W^{2,q}(\mathbb{T})\times \mathbb{R})_{1-1/p,p}\\
    &=(L^q(\mathbb{T}),W^{2,q}(\mathbb{T}))_{1-1/p,p}\times (\mathbb{R},\mathbb{R})_{1-1/p,p}\\
    &=B^{2-2/p}_{q,p}(\mathbb{T})\times \mathbb{R}
\end{align*}
its norm is given by 
\[
\|(f,L)\|_{X_p}=\|f\|\p+|L|
\]
for all $(f,L)\in X_p.$

We choose $p,q$ such that $1>2/p+1/q$,  by the embedding theorem in \cite[Section 2.8.1]{Tri78}, the continuous embedding $B^{2-2/p}_{q,p}(\mathbb{T})\hookrightarrow C^{1,\alpha}(\mathbb{T})$ holds for some $\alpha>0$, i.e. there exists a constant $C>0$ such that $\|f\|_{C^{1,\alpha}(\mathbb{T})}\leq C\|f\|\p$ for all $f\in B^{2-2/p}_{q,p}(\mathbb{T})$, thus, we have 
\[
\|f\|_{\infty}, \|\partial_rf\|_{\infty}\leq C\|f\|_{B^{2-2/p}_{q,p}(\mathbb{T})}
\]

To construct a solution of (\ref{quasilinear equation}), we investigate the following truncated equation
 \begin{equation}\label{truncated equation}
    \left\{\begin{aligned}
        d(f(t),L(t))^{\top}&=\left[-A_n(f(t),L(t))(f(t),L(t))^{\top}+F_n(t,f(t),L(t))\right]dt+B_n(f(t),L(t))dW_t\\
        (f(0),L(0))&=(k_0(rL_0),L_0)
    \end{aligned}\right.
\end{equation}

where $A_n(g,M)=A(g,T_nM)$, $F_n(g,M)=F(g,T_nM)$, $G_n(g,M)=G(g,T_nM)$. The cut-off mapping $T_n$ are defined by 
\begin{equation}
    \quad T_nM:=\left\{\begin{array}{ll}
        \dfrac{M}{n|M|}, & 0<|M|<\dfrac{1}{n},\\
        M, &  |M|\geq \dfrac{1}{n} .
    \end{array}\right.
\end{equation}
The idea is to use such a truncation to extend global nonlinearity to local ones, and it was used several times to solve semilinear equations(see \cite{Seidler1993,vNVW08,Brz95,BMS05}).

\begin{lemma}\label{lemma}
    Given $n\in \mathbb{N}$, the above cut-off mapping are of linear growth and Lipschitz continuous, i.e.
    \[
    \|(g,T_nM)\|_{X_p}\leq 1+\|(g,M)\|_{X_p}.
    \]\[
\|(g_1,T_nM_1)-(g_2,T_nM_2)\|_{X_p}\leq \|(g_1,M_1)-(g_2,M_2)\|_{X_p}.
\]
\end{lemma}

\begin{proof}[Proof of Lemma \ref{lemma}]
Since for all $n\in\mathbb{N}$,
\begin{equation*}
    |T_nM|\leq \left\{\begin{array}{cc}
        |M| & \text{if } |M|>1\\
        1 & \text{if } |M|\leq 1
    \end{array}\right.
\end{equation*}
we can get 
\[
\|(R_ng,T_nM)\|_{X_p}\leq 1+\|(g,M)\|_{X_p}
\]

It is easy to show that 
\[
|T_nM_1-T_nM_2|\leq |M_1-M_2|
\]
Thus, we have 
\[\|(g_1,T_nM_1)-(g_2,T_nM_2)\|_{X_p}\leq \|(g_1,M_1)-(g_2,M_2)\|_{X_p}.
\]
\end{proof}

\begin{lemma}
    For initial data $(f_0,L_0)\in B^{2-2/p}_{q,p}(\mathbb{T})\times \R_+ \subset X_p$, the initial operator
    \begin{equation}
    A(f_0,L_0)=-\left[\left(2\sigma^2\pi^2r^2+\frac{1}{L_0^2}\right)\2\right]I_2\in\mathcal{SMR}^{\bullet}_p(T)\footnote{for the definition of $\mathcal{SMR}^{\bullet}_p(T)$, we refer to \cite[Definition 3.5]{AV22a}}.
\end{equation}
\end{lemma}
\begin{proof}
    First, we choose $\mu\in (0,\frac{\pi}{2}]$.
    Since $L_0$ denotes the length of the initial closed curve $\gamma_0$, and $f_0(r)=k_0(rL_0)$ denotes the curvature of $\gamma_0$, both of them must be bounded, i.e. there existes a constant $M>0$, such that 
    \[
    \frac{1}{M}\leq L_0\leq M,\quad \sup_{\T}|f_0|\leq M.
    \]

    Note that the initial operator $A(f_0,L_0)$ can be expressed in the following way,   
    \[
    A(f_0,L_0)=a_2(r)D^2.
    \]
    where $D=-i\1$. Then, the coefficient of $A(f_0,L_0)$
    \[
\|a_2\|\f=\left\|\left(2\sigma^2\pi^2r^2+\frac{1}{L_0^2}\right)I_2\right\|\f\leq 2\sigma^2\pi^2+M^2.
\]
It is clear that $a_2\in \mathrm{BUC}(\mathbb{T},L(\R^2))$. 
Then the principal symbol of $A(f_0,L_0)$
\[
(A(f_0,L_0))_{\pi}(\xi)=\left(2\sigma^2\pi^2r^2+\frac{1}{L_0^2}\right)I_2>0,
\]
Its spectrum satisfies
\[\sigma((A(f_0,L_0))_{\pi}(\xi))=\left\{2\sigma^2\pi^2r^2+\frac{1}{L_0^2}I_2\right\}\subset \Sigma_{\mu}.\]
Its inverse satisfies
\[
\left[(A(f_0,L_0))_{\pi}(\xi)\right]^{-1}=\left(2\sigma^2\pi^2r^2+\frac{1}{L_0^2}\right)^{-1}I_2.
\]
Hence, 
\[
\det\left[A(f_0,L_0)_{\pi}(\xi)\right]^{-1}=\left(2\sigma^2\pi^2r^2+\frac{1}{L_0^2}\right)^{-2}\leq L^4_0\leq  M^4.
\]
Thus,   $A(f_0,L_0)$ is a uniformly elliptic operator.  
By \cite[Theorem 6.1]{DS97}, there exists a constant $s>0$ such that operator
\[
sI+A(f_0,L_0)
\]
has a bounded $H^{\infty}$-calculus of angle less than $\pi/2$. By Theorem 3.7 in \cite{AV22a}, we have 
\[
A(f_0,L_0)\in \mathcal{SMR}^{\bullet}_p(T).
\]   
\end{proof}

\begin{lemma}
    $A_n, F_n$ and $B_n$ satisfy (HA), (HF) and (HG) in Section 4.1 in \cite{AV22a} respectively.
\end{lemma}
\begin{proof}
For all $n\in\mathbb{N}$, assume that $(f,L),(f_i,L_i)\in B^{2-2/p}_{q,p}(\mathbb{T})\times \mathbb{R}, i=1,2$ satisfy 
\begin{align*}
    &\|f\|\q+|L|=\|(f,L)\|_{X_p}\leq n,\\
    &\|f_i\|\q+|L_i|=\|(f_i,L_i)\|_{X_p}\leq n,\\
    &|L|,|L_i|\geq1/n, \quad i=1,2.
\end{align*}
For any $(g,M)\in W^{2,q}(\T)\times \R$,
\begin{align*}
    \|A(f,L)(g,M)^{\top}\|_{X}&=\left\|-2\sigma^2\pi^2 r^2\2g-\frac{1}{L^2}\2 g\right\|\p\\
    &\leq (2\sigma^2\pi^2+n^2)\|\2g\|\p\\
    &\leq C_A(n)\|(g,M)\|_{X_1}.
\end{align*}

\begin{align*}
    \left\|(A(f_1,L_1)-A(f_2,L_2)(g,M)^{\top})\right\|_X&=\left\|\left(\frac{1}{L_2^2}-\frac{1}{L^2_1}\right)\2g\right\|\p\\
    &\leq \frac{|L_1+L_2||L_1-L_2|}{L^2_1L^2_2}\|\2g\|\p\\
    &\leq 2n^3|L_1-L_2|\|(g,M)\|_{X}\\
    &\leq L_A(n)\|(f_1,L_1)-(f_2,L_2)\|_{X_p}\|(g,M)\|_{X}.
\end{align*}
Operator $A(f,L)$  is Lipschitz continuous for $(f,L)$ satisfying $\|f\|\p\leq n$ and  $1/n\leq|L|\leq n$.

As for the nonlinearity $F$: 
\begin{align*}
    &\|F(t,f,L)\|_{X}\\
    &=\left\|-r(\partial_s f)\int_{\mathbb{T}} f^2dr+f^3+\sigma^2\left(-4\pi rLf\1f +2\pi^2r\1f+f^3L^2-\pi f^2L\right) \right\|\p\\
    &\quad +\left|2\sigma^2\pi^2L-L\int_{\mathbb{T}} f^2dr\right|\\
    &\leq \|\partial_s f\|_{\infty}\|f\|_{\infty}^2 + \|f\|_{\infty}^3 +4\pi \sigma^2|L|\|f\|\f\|\1f\|\f+2\sigma^2\pi^2\|\1f\|\f  \\
    &\quad +\sigma^2\|f\|^3\f|L|^2++\sigma^2\pi \|f\|^2\f|L|+2\sigma^2\pi^2|L|+|L|\|f\|^2\f \\
    &\leq (2+\sigma^2n^2)C^3\|f\|_{B^{2-2/p}_{q,p}}^3+(5\sigma^2\pi+1)C^2 \|f\|_{B^{2-2/p}_{q,p}}^2|L|+2\sigma^2\pi^2C\|f\|_{B^{2-2/p}_{q,p}}+2\sigma^2\pi^2|L|\\
    &\leq  \left[C^3(\sigma^2n^4+n^2)+2C\sigma^2\pi^2\right]\|f\|_{B^{2-2/p}_{q,p}}+\left[C^2(5\sigma^2\pi n^2+n^2)+2\sigma^2\pi^2\right]|L| \\
    &\leq C_F(n)\left(\|f\|_{B^{2-2/p}_{q,p}}+|L|\right) \\
    &= C_F(n)\|(f,L)\|_{X_p}.
\end{align*}
and
\begin{align*}
    \|F(t,&f_1,L_1)-F(t,f_2,L_2)\|_X \\
    &= \left|2\sigma^2\pi^2(L_1-L_2)-L_1\int_{\mathbb{T}} f_1^2dr+L_2\int_{\mathbb{T}} f_2^2dr\right| \\
    &=\left\|-r\1f\int_{\T}f^2dr+f^3+\sigma^2\left(-4\pi rLf\1f+2\pi^2 r \1f +f^3L^2-\pi f^2L\right)\right.  \\
    &\quad +\left.\left(-r\1f\int_{\T}f^2dr+f^3+\sigma^2\left(-4\pi rLf\1f+2\pi^2 r \1f +f^3L^2-\pi f^2L\right)\right)\right\|_{\p}\\
    &\leq 2\sigma^2\pi^2|L_1-L_2| + |L_1-L_2|\|f_1\|_{L^\infty}^2 + |L_2|\|f_1+f_2\|_{L^\infty}\|f_1-f_2\|_{L^\infty} \\
    &\quad + \|\partial_s f_1-\partial_s f_2\|_{L^\infty}\|f_1\|_{L^\infty}^2 + \|\partial_s f_2\|_{L^\infty}\|f_1+f_2\|_{L^\infty}\|f_1-f_2\|_{L^\infty} \\
    &\quad + (1+\sigma^2L_1^2)(\|f_1\|_{L^\infty}^2 + \|f_1\|_{L^\infty}\|f_2\|_{L^\infty} + \|f_2\|_{L^\infty}^2)\|f_1-f_2\|_{L^\infty} \\
    &\quad + \sigma^2\|f_2\|_{L^\infty}^3|L_1+L_2||L_1-L_2| + \sigma^2\pi \|f_1+f_2\|_{L^\infty}\|f_1-f_2\|_{L^\infty}|L_1| \\
    &\quad + \sigma^2\pi\|f_2\|_{L^\infty}^2|L_1-L_2| +2\sigma \pi^2\|f_1-f_2\|\f+4\sigma^2\pi\|f_1\|\f\|\1f\|\f|L_1-L_2|\\
    &\quad ++4\sigma^2\pi |L_2|\|\1f_1-\1f_2\|\f\|f_1\|\f+4\sigma^2\pi |L_2|\|\2f\|\f\|f_1-f_2\|\f\\
    &\leq \left(2\sigma^2\pi^2+C^2n^2+2C^3\sigma^2n^4+5C^2\sigma^2\pi\right)|L_1-L_2|\\
    &\quad + \left(2C^2n^2+6C^3n^2+3C^3\sigma^2n^4+2C\sigma^2\pi n^2+2C\sigma\pi^2\right)\|f_1-f_2\|_{B^{2-2/p}_{q,p}}  \\
    &\leq L_F(n)\left(\|f_1-f_2\|_{B^{2-2/p}_{q,p}}+|L_1-L_2|\right) \\
    &= L_F(n)\|(f_1,L_1)-(f_2,L_2)\|_{X_p}.
\end{align*}

Note that $X^{1/2}=W^{1,q}(\mathbb{T})\times \mathbb{R}$, by the definition of the space $\gamma(\mathbb{R};W^{1,q}(\mathbb{T}))$, it is easy to prove that for any $g\in W^{1,q}(\mathbb{T})$ we have $\|g\|_{{\gamma(\mathbb{R}}; W^{1,q}(\mathbb{T}))}=\|g\|_{{W^{1,q}(\mathbb{T})}}$.

\begin{align*}
    \|B_1(f,L)\|_{\gamma(\mathbb{R};X^{1/2})} &= \left\|\sigma f^2L\right\|_{W^{1,q}(\mathbb{T})} \\
    &= \sigma|L|\left( \|f^2\|_{L^q} + \|2f \partial_s f\|_{L^q} \right) \\
    &\leq \sigma |L| \|f\|_{L^\infty}^2 + 2\sigma |L| \|f\|_{L^\infty} \|\partial_s f\|_{L^\infty} \\
    &\leq \sigma C^2 n^2 |L| + 2\sigma C^2 n^2 \|f\|_{B^{2-2/p}_{q,p}} \\
    &\leq 2\sigma C^2 n^2 \left( \|f\|_{B^{2-2/p}_{q,p}} + |L| \right) \\
    &= C_{B_1}(n) \|(f,L)\|_{X_p}.
\end{align*}
and
\begin{align*}
    \|B_1(f_1,&L_1)-B_1(f_2,L_2)\|_{\gamma(\mathbb{R};X^{1/2})} \\
    &= \sigma \| f_1^2 L_1 - f_2^2 L_2 \|_{W^{1,q}(\mathbb{T})} \\
    &= \sigma \| f_1^2 L_1 - f_2^2 L_2 \|_{L^q} + \sigma \| 2f_1 (\partial_s f_1) L_1 - 2f_2 (\partial_s f_2) L_2 \|_{L^q} \\
    &\leq \sigma \|f_1 + f_2\|_{L^\infty} \|f_1 - f_2\|_{L^\infty} |L_1| + \sigma \|f_2\|_{L^\infty}^2 |L_1 - L_2| \\
    &\quad + 2\sigma \|f_1 - f_2\|_{L^\infty} \|\partial_s f_1\|_{L^\infty} |L_1| + 2\sigma \|f_2\|_{L^\infty} \|\partial_s f_1 - \partial_s f_2\|_{L^\infty} |L_1| \\
    &\quad + 2\sigma \|f_2\|_{L^\infty} \|\partial_s f_2\|_{L^\infty} |L_1 - L_2| \\
    &\leq 6\sigma C^2 n^2 \|f_1 - f_2\|_{B^{2-2/p}_{q,p}} + 3\sigma C^2 n^2 |L_1 - L_2| \\
    &\leq 6\sigma C^2 n^2 \left( \|f_1 - f_2\|_{B^{2-2/p}_{q,p}} + |L_1 - L_2| \right) \\
    &= L_{B_1}(n) \|(f_1,L_1) - (f_2,L_2)\|_{X_p}.
\end{align*}
and 
\begin{align*}
    \|B_1(f_1,&L_1)-B_1(f_2,L_2)\|_{\gamma(\mathbb{R};X^{1/2})} \\
    &= \sigma \| f_1^2 L_1 - f_2^2 L_2 \|_{W^{1,q}(\mathbb{T})} \\
    &= \sigma \| f_1^2 L_1 - f_2^2 L_2 \|_{L^q} + \sigma \| 2f_1 (\partial_s f_1) L_1 - 2f_2 (\partial_s f_2) L_2 \|_{L^q} \\
    &\leq \sigma \|f_1 + f_2\|_{L^\infty} \|f_1 - f_2\|_{L^\infty} |L_1| + \sigma \|f_2\|_{L^\infty}^2 |L_1 - L_2| \\
    &\quad + 2\sigma \|f_1 - f_2\|_{L^\infty} \|\partial_s f_1\|_{L^\infty} |L_1| + 2\sigma \|f_2\|_{L^\infty} \|\partial_s f_1 - \partial_s f_2\|_{L^\infty} |L_1| \\
    &\quad + 2\sigma \|f_2\|_{L^\infty} \|\partial_s f_2\|_{L^\infty} |L_1 - L_2| \\
    &\leq 6\sigma C^2 n^2 \|f_1 - f_2\|_{B^{2-2/p}_{q,p}} + 3\sigma C^2 n^2 |L_1 - L_2| \\
    &\leq 6\sigma C^2 n^2 \left( \|f_1 - f_2\|_{B^{2-2/p}_{q,p}} + |L_1 - L_2| \right) \\
    &= L_{B_1}(n) \|(f_1,L_1) - (f_2,L_2)\|_{X_p}.
\end{align*}
Besides, for all $(f_1,L_1),(f_2,L_2),(f,L)\in X_1$, we have
\begin{align*}
    \|B_2(f,L)\|_{\gamma(\mathbb{R};X^{1/2})} 
    &= \| (-2\sigma\pi r \partial_s f, -2\sigma\pi L) \|_{\gamma(\mathbb{R};W^{1,q}(\mathbb{T})\times\mathbb{R})} \\
    &= 2\sigma\pi \left( \| r \partial_s f \|_{\gamma(\mathbb{R};W^{1,q}(\mathbb{T})} + |L|_{\gamma(\mathbb{R};\mathbb{R})} \right) \\
    &= 2\sigma\pi \left( \| r \partial_s f \|_{W^{1,q}(\mathbb{T})} + |L| \right) \\
    &= 2\sigma\pi \left( \| r \partial_s f \|_{L^q} + \| \partial_s f + r \partial_s^2 f \|_{L^q} + |L| \right) \\
    &\leq 2\sigma\pi \left( 2 \| \partial_s f \|_{L^q} + \| \partial_s^2 f \|_{L^q} + |L| \right) \\
    &\leq 4\sigma\pi \left( \|f\|_{W^{2,q}(\mathbb{T})} + |L| \right) \\
    &= 4\sigma\pi \|(f,L)\|_{X_1}.\\
    \\
    \|B_2(f_1,&L_1) - B_2(f_2,L_2)\|_{\gamma(\mathbb{R};X^{1/2})} \\
    &= 2\sigma\pi \| \left( r \partial_s (f_1 - f_2), L_1 - L_2 \right) \|_{\gamma(\mathbb{R};W^{1,q}(\mathbb{T})\times\mathbb{R})} \\
    &\leq 4\sigma\pi \|(f_1,L_1) - (f_2,L_2)\|_{X_1}.
\end{align*}

\end{proof}

We apply Theorem 4.5 in \cite{AV22a} to the truncated equation (\ref{truncated equation}) and obtain for every $n\in\mathbb{N}$ a $L^p$-unique maximal local solution $((f_n,L_n),\tau_n)$.  And $(f_n,L_n)$ satisfies \eqref{quasilinear equation} on $[0,\sigma_n)$, where 
 \begin{equation}
    \sigma_n:=\tau_n\wedge\inf\left\{t\in[0,\tau_n):\|f_n\|\q>n\text{\, or \,} L\in (0,1/n)\cup(n,\infty)\right\}
\end{equation}
It is easy to see that $\sigma_n$ is a $\mathcal{F}$-stopping time. By Lemma 4.10 in \cite{Hornung}, we have : for $\mathbb{P}-a.s. \omega\in \Omega$, $(\sigma_n(\omega))_n$ is monotonously increasing beginning from some $n=n(\omega)\in \mathbb{N}$. Moreover, for all $l>k\geq n(\omega)$  and $t\in [0,\sigma_n(\omega))$ we have 
\[
(f_k,L_k)(\omega,t)=(f_l,L_l)(\omega,t).
\]
Then,  applying the same method for proving Theorem 4.11 in \cite{Hornung}, we can get $L^p$-maximal local solution for \eqref{quasilinear equation} and the corresponding bolw-up criterion ( or we can apply directly Theorem 4.9 in \cite{AV22b}). Finally, we get the following theorem :
\begin{theorem}
    For given $p,q>0$ satisfying $1>2/p+1/q$, there exist a constant $\epsilon>0$ such that  for all $\sigma$ with $0<\sigma\leq \epsilon$ The following equations
    \begin{equation}
        \left\{
    \begin{aligned}
        df=&\left[\left(2\sigma^2\pi^2r^2+\frac{1}{L^2}\right)\2f-4\sigma^2\pi rLf\1f+2\sigma^2\pi^2 r \1f -r\1f\int_{\T} f^2dr+f^3+\sigma^2f^3L^2-\sigma^2\pi f^2L\right]dt\\
        &\quad +\sigma\left(f^2L-2\pi r\1f\right)dW_t\\
        dL=&L\left(2\sigma^2\pi^2-\int_{\T} f^2dr\right)dt -2\sigma\pi LdW_t\\
        f(r,0)&=k_0(rL_0),\quad L(0)=L_0, \quad r\in\mathbb{T}
    \end{aligned}
    \right.
    \end{equation}
    have a unique $L^p$-maximal local solution $(f,L,\tau)$ such that $\tau>0 $ a.s.. For $(f,L,\tau)$, there exists a sequence of stopping time $(\tau_n)_{n\geq 1}$ such that $\tau>0$ a.s., and for all $n\geq 1, \theta\in [0,1/2)$, we have 
\[
(f,L)\in L^p(\Omega;H^{\theta,p}([0,\tau_n];H^{2(1-\theta),q}(\T)\times \R))\cap L^p(\Omega;C([0,\tau_n]; B^{2-2/p}_{q,p}(\T)\times\R)).
\]

Besides, the following blow-up criterion holds.
\begin{equation*}
    \mathbb{P} \left\{\begin{aligned}
         &\tau<T,  \,\|f\|_{L^p(0,\tau;H^{2-2\theta,q}(\T))}<\infty,\, 0<L<\infty, \,\\
         &(f,L):[0,\tau)\to B^{2-2/p}_{q,p}(\T)\times\mathbb{R}_+\,\text{ is uniformly continuous }
     \end{aligned}\right\}=0.
 \end{equation*}
    
\end{theorem}

Then, we get the theorem for $(k,L)$:
\begin{theorem}
    For given $p,q>0$ satisfying $1>2/p+1/q$, there exist a constant $\epsilon>0$ such that  for all $\sigma$ with $0<\sigma\leq \epsilon$, the stochastic curve shortening flow driven by scale-dependent noise 
    \begin{equation}
    \left\{
    \begin{aligned}
        & dk(t)=\left(\partial_{ss}k+k^3\right)dt+\sigma k^2L \circ dW_t,\quad\quad s\in [0,L(t)],\\
       &dL(t)=-\int^{L(t)}_0k^2dsdt-2\sigma\pi L\circ dW_t,\\
       &k(s,0)=k_0(s), \quad L(0)=L_0, \quad s\in [0,L_0].
    \end{aligned}\right.
\end{equation}
have a  unique $L^p$-maximal local solution $(k,L,\tau)$ such that $\tau>0 $ a.s.. For $(k,L,\tau)$, there exists a sequence of stopping time $(\tau_n)_{n\geq 1}$ such that $\tau>0$ a.s., and for all $n\geq 1, \theta\in [0,1/2)$, we have 
\[
(k,L)\in L^p(\Omega;H^{\theta,p}([0,\tau_n];H^{2(1-\theta),q}([0,L])\times \R))\cap L^p(\Omega;C([0,\tau_n]; B^{2-2/p}_{q,p}([0,L])\times\R)).
\]

Besides, the following blow-up criterion holds.
\begin{equation*}
    \mathbb{P} \left\{\begin{aligned}
         &\tau<T,  \,\|f\|_{L^p(0,\tau;H^{2-2\theta,q}([0,L]))}<\infty,\, 0<L<\infty, \,\\
         &(f,L):[0,\tau)\to B^{2-2/p}_{q,p}([0,L])\times\mathbb{R}_+\,\text{ is uniformly continuous }
     \end{aligned}\right\}=0.
 \end{equation*}
\end{theorem}

 \begin{remark}
     By the construction of the maximal stopping time $\tau$ in \cite{Hornung}, the blow-up time corresponds to the time when $k\to\infty$ or $L\to 0$, it means that the flow develops singularities or shrinks to a point. 
 \end{remark}
 \begin{remark}
     By the expression 
     \[
     L(t)=L_0\exp\left(-\int^t_0\int_{\T} f^2(r,\tau)drd\tau-2\sigma\pi W_t\right)
     \]
     We can get that 
     \[
     \left\{L(t)=0\right\}=\left\{\sup_{\mathbb{T}}\|f(t)\|=\infty\right\}\bigcup\left\{W_t=+\infty\right\}=\left\{\sup\|k(t)\|=\infty\right\}\cup\{W_t=+\infty\}
     \]
     \[
     \left\{L(t)=+\infty\right\}=\left\{W_t=-\infty\right\}
     \]
 \end{remark}

\bibliographystyle{alpha}
\bibliography{sample}

@article {Hornung,
    AUTHOR = {Hornung, Luca},
     TITLE = {Quasilinear parabolic stochastic evolution equations via
              maximal {$L^p$}-regularity},
   JOURNAL = {Potential Anal.},
  FJOURNAL = {Potential Analysis. An International Journal Devoted to the
              Interactions between Potential Theory, Probability Theory,
              Geometry and Functional Analysis},
    VOLUME = {50},
      YEAR = {2019},
    NUMBER = {2},
     PAGES = {279--326},
      ISSN = {0926-2601,1572-929X},
   MRCLASS = {60H15 (35B65 35K59 35R60 58D25 60H30 65J08 76A05)},
  MRNUMBER = {3905531},
       DOI = {10.1007/s11118-018-9683-9},
       URL = {https://doi.org/10.1007/s11118-018-9683-9},
}

@book {DJT,
    AUTHOR = {Diestel, Joe and Jarchow, Hans and Tonge, Andrew},
     TITLE = {Absolutely summing operators},
    SERIES = {Cambridge Studies in Advanced Mathematics},
    VOLUME = {43},
 PUBLISHER = {Cambridge University Press, Cambridge},
      YEAR = {1995},
     PAGES = {xvi+474},
      ISBN = {0-521-43168-9},
   MRCLASS = {46-02 (46Bxx 46M05 47-02 47B10 47D50)},
  MRNUMBER = {1342297},
MRREVIEWER = {Andreas\ Defant},
       DOI = {10.1017/CBO9780511526138},
       URL = {https://doi.org/10.1017/CBO9780511526138},
}

@incollection {vanN,
    AUTHOR = {van Neerven, Jan},
     TITLE = {{$\gamma$}-radonifying operators---a survey},
 BOOKTITLE = {The {AMSI}-{ANU} {W}orkshop on {S}pectral {T}heory and
              {H}armonic {A}nalysis},
    SERIES = {Proc. Centre Math. Appl. Austral. Nat. Univ.},
    VOLUME = {44},
     PAGES = {1--61},
 PUBLISHER = {Austral. Nat. Univ., Canberra},
      YEAR = {2010},
      ISBN = {0-7315-5208-3},
   MRCLASS = {47B10 (28C20 46B09 60B11 60H05)},
  MRNUMBER = {2655391},
MRREVIEWER = {Werner\ Linde},
}

@article {vNVW07,
    AUTHOR = {van Neerven, J. M. A. M. and Veraar, M. C. and Weis, L.},
     TITLE = {Stochastic integration in {UMD} {B}anach spaces},
   JOURNAL = {Ann. Probab.},
  FJOURNAL = {The Annals of Probability},
    VOLUME = {35},
      YEAR = {2007},
    NUMBER = {4},
     PAGES = {1438--1478},
      ISSN = {0091-1798,2168-894X},
   MRCLASS = {60H05 (28C20 60B11)},
  MRNUMBER = {2330977},
MRREVIEWER = {Martin\ Ondrej\'at},
       DOI = {10.1214/009117906000001006},
       URL = {https://doi.org/10.1214/009117906000001006},
}

@article {DHP,
    AUTHOR = {Denk, Robert and Hieber, Matthias and Pr\"uss, Jan},
     TITLE = {{$R$}-boundedness, {F}ourier multipliers and problems of
              elliptic and parabolic type},
   JOURNAL = {Mem. Amer. Math. Soc.},
  FJOURNAL = {Memoirs of the American Mathematical Society},
    VOLUME = {166},
      YEAR = {2003},
    NUMBER = {788},
     PAGES = {viii+114},
      ISSN = {0065-9266,1947-6221},
   MRCLASS = {35-02 (35J30 35K25 42B15)},
  MRNUMBER = {2006641},
MRREVIEWER = {Alain\ Brillard},
       DOI = {10.1090/memo/0788},
       URL = {https://doi.org/10.1090/memo/0788},
}

@article {CPS,
    AUTHOR = {Cl\'ement, P. and de Pagter, B. and Sukochev, F. A. and
              Witvliet, H.},
     TITLE = {Schauder decomposition and multiplier theorems},
   JOURNAL = {Studia Math.},
  FJOURNAL = {Studia Mathematica},
    VOLUME = {138},
      YEAR = {2000},
    NUMBER = {2},
     PAGES = {135--163},
      ISSN = {0039-3223,1730-6337},
   MRCLASS = {47A65 (42A45 46E99)},
  MRNUMBER = {1749077},
MRREVIEWER = {Maria\ J.\ Carro},
}

@incollection {KW,
    AUTHOR = {Kunstmann, Peer C. and Weis, Lutz},
     TITLE = {Maximal {$L_p$}-regularity for parabolic equations, {F}ourier
              multiplier theorems and {$H^\infty$}-functional calculus},
 BOOKTITLE = {Functional analytic methods for evolution equations},
    SERIES = {Lecture Notes in Math.},
    VOLUME = {1855},
     PAGES = {65--311},
 PUBLISHER = {Springer, Berlin},
      YEAR = {2004},
      ISBN = {3-540-23030-0},
   MRCLASS = {47D06 (34G10 35D10 35J55 35K20 35K90 42B20 47A60)},
  MRNUMBER = {2108959},
MRREVIEWER = {Xuan\ Thinh\ Duong},
       DOI = {10.1007/978-3-540-44653-8\_2},
       URL = {https://doi.org/10.1007/978-3-540-44653-8_2},
}

@article {GH86,
    AUTHOR = {Gage, M. and Hamilton, R. S.},
     TITLE = {The heat equation shrinking convex plane curves},
   JOURNAL = {J. Differential Geom.},
  FJOURNAL = {Journal of Differential Geometry},
    VOLUME = {23},
      YEAR = {1986},
    NUMBER = {1},
     PAGES = {69--96},
      ISSN = {0022-040X,1945-743X},
   MRCLASS = {53A04 (35K05 52A40 58E99 58G11)},
  MRNUMBER = {840401},
MRREVIEWER = {R.\ Osserman},
       URL = {http://projecteuclid.org/euclid.jdg/1214439902},
}

@article {DS97,
    AUTHOR = {Duong, Xuan T. and Simonett, Gieri},
     TITLE = {{$H_\infty$}-calculus for elliptic operators with nonsmooth
              coefficients},
   JOURNAL = {Differential Integral Equations},
  FJOURNAL = {Differential and Integral Equations. An International Journal
              for Theory \& Applications},
    VOLUME = {10},
      YEAR = {1997},
    NUMBER = {2},
     PAGES = {201--217},
      ISSN = {0893-4983},
   MRCLASS = {47F05 (35J45 47A60 47N20)},
  MRNUMBER = {1424807},
MRREVIEWER = {Vladimir\ B.\ Vasilyev},
}

@article{Seidler1993,
author = {Seidler, Jan},
journal = {Mathematica Bohemica},
language = {eng},
number = {1},
pages = {67-106},
publisher = {Institute of Mathematics, Academy of Sciences of the Czech Republic},
title = {Da Prato-Zabczyk's maximal inequality revisited. I.},
url = {http://eudml.org/doc/29167},
volume = {118},
year = {1993},
}

@article {vNVW08,
    AUTHOR = {van Neerven, J. M. A. M. and Veraar, M. C. and Weis, L.},
     TITLE = {Stochastic evolution equations in {UMD} {B}anach spaces},
   JOURNAL = {J. Funct. Anal.},
  FJOURNAL = {Journal of Functional Analysis},
    VOLUME = {255},
      YEAR = {2008},
    NUMBER = {4},
     PAGES = {940--993},
      ISSN = {0022-1236,1096-0783},
   MRCLASS = {35R60 (35B65 35K55 46N30 47D06 60H15 60H20)},
  MRNUMBER = {2433958},
MRREVIEWER = {Krystyna\ Twardowska},
       DOI = {10.1016/j.jfa.2008.03.015},
       URL = {https://doi.org/10.1016/j.jfa.2008.03.015},
}

@article {Brz95,
    AUTHOR = {Brze\'zniak, Zdzis\l aw},
     TITLE = {Stochastic partial differential equations in {M}-type {$2$}
              {B}anach spaces},
   JOURNAL = {Potential Anal.},
  FJOURNAL = {Potential Analysis. An International Journal Devoted to the
              Interactions between Potential Theory, Probability Theory,
              Geometry and Functional Analysis},
    VOLUME = {4},
      YEAR = {1995},
    NUMBER = {1},
     PAGES = {1--45},
      ISSN = {0926-2601,1572-929X},
   MRCLASS = {35R60 (34F05 34G10 35R20 47N20 60H15)},
  MRNUMBER = {1313905},
MRREVIEWER = {Stanis\l aw\ W\polhk edrychowicz},
       DOI = {10.1007/BF01048965},
       URL = {https://doi.org/10.1007/BF01048965},
}

@article {BMS05,
    AUTHOR = {Brze\'zniak, Zdzis\l aw and Maslowski, Bohdan and Seidler,
              Jan},
     TITLE = {Stochastic nonlinear beam equations},
   JOURNAL = {Probab. Theory Related Fields},
  FJOURNAL = {Probability Theory and Related Fields},
    VOLUME = {132},
      YEAR = {2005},
    NUMBER = {1},
     PAGES = {119--149},
      ISSN = {0178-8051,1432-2064},
   MRCLASS = {60H15 (35Q72 35R60 74H50 74K10)},
  MRNUMBER = {2136869},
MRREVIEWER = {Grigorios\ A.\ Pavliotis},
       DOI = {10.1007/s00440-004-0392-5},
       URL = {https://doi.org/10.1007/s00440-004-0392-5},
}

@article {EvR12,
    AUTHOR = {Es-Sarhir, Abdelhadi and von Renesse, Max-K.},
     TITLE = {Ergodicity of stochastic curve shortening flow in the plane},
   JOURNAL = {SIAM J. Math. Anal.},
  FJOURNAL = {SIAM Journal on Mathematical Analysis},
    VOLUME = {44},
      YEAR = {2012},
    NUMBER = {1},
     PAGES = {224--244},
      ISSN = {0036-1410,1095-7154},
   MRCLASS = {47D07 (35R60 60H15)},
  MRNUMBER = {2888287},
MRREVIEWER = {Markus\ Kunze},
       DOI = {10.1137/100798235},
       URL = {https://doi.org/10.1137/100798235},
}

@article {HRvR17,
    AUTHOR = {Hofmanov\'a, Martina and R\"oger, Matthias and von Renesse,
              Max},
     TITLE = {Weak solutions for a stochastic mean curvature flow of
              two-dimensional graphs},
   JOURNAL = {Probab. Theory Related Fields},
  FJOURNAL = {Probability Theory and Related Fields},
    VOLUME = {168},
      YEAR = {2017},
    NUMBER = {1-2},
     PAGES = {373--408},
      ISSN = {0178-8051,1432-2064},
   MRCLASS = {60H15 (53C44)},
  MRNUMBER = {3651056},
       DOI = {10.1007/s00440-016-0713-5},
       URL = {https://doi.org/10.1007/s00440-016-0713-5},
}

@article {DHR21,
    AUTHOR = {Dabrock, Nils and Hofmanov\'a, Martina and R\"oger, Matthias},
     TITLE = {Existence of martingale solutions and large-time behavior for
              a stochastic mean curvature flow of graphs},
   JOURNAL = {Probab. Theory Related Fields},
  FJOURNAL = {Probability Theory and Related Fields},
    VOLUME = {179},
      YEAR = {2021},
    NUMBER = {1-2},
     PAGES = {407--449},
      ISSN = {0178-8051,1432-2064},
   MRCLASS = {60H15 (53E10 60H30)},
  MRNUMBER = {4221662},
       DOI = {10.1007/s00440-020-01012-6},
       URL = {https://doi.org/10.1007/s00440-020-01012-6},
}

@book {Tri78,
    AUTHOR = {Triebel, H.},
     TITLE = {Interpolation theory, function spaces, differential operators},
 PUBLISHER = {VEB Deutscher Verlag der Wissenschaften, Berlin},
      YEAR = {1978},
     PAGES = {528},
   MRCLASS = {46E35 (35Jxx 46M35)},
  MRNUMBER = {500580},
MRREVIEWER = {Robert\ D.\ Brown},
}

@article {Hui84,
    AUTHOR = {Huisken, Gerhard},
     TITLE = {Flow by mean curvature of convex surfaces into spheres},
   JOURNAL = {J. Differential Geom.},
  FJOURNAL = {Journal of Differential Geometry},
    VOLUME = {20},
      YEAR = {1984},
    NUMBER = {1},
     PAGES = {237--266},
      ISSN = {0022-040X,1945-743X},
   MRCLASS = {53C45 (49F05 58F17)},
  MRNUMBER = {772132},
MRREVIEWER = {R.\ Schneider},
       URL = {http://projecteuclid.org/euclid.jdg/1214438998},
}

@article {Hui90,
    AUTHOR = {Huisken, Gerhard},
     TITLE = {Asymptotic behavior for singularities of the mean curvature
              flow},
   JOURNAL = {J. Differential Geom.},
  FJOURNAL = {Journal of Differential Geometry},
    VOLUME = {31},
      YEAR = {1990},
    NUMBER = {1},
     PAGES = {285--299},
      ISSN = {0022-040X,1945-743X},
   MRCLASS = {53A10 (35B99 53C45 58G11)},
  MRNUMBER = {1030675},
MRREVIEWER = {Dennis\ M.\ DeTurck},
       URL = {http://projecteuclid.org/euclid.jdg/1214444099},
}

@article {CM12,
    AUTHOR = {Colding, Tobias H. and Minicozzi, II, William P.},
     TITLE = {Generic mean curvature flow {I}: generic singularities},
   JOURNAL = {Ann. of Math. (2)},
  FJOURNAL = {Annals of Mathematics. Second Series},
    VOLUME = {175},
      YEAR = {2012},
    NUMBER = {2},
     PAGES = {755--833},
      ISSN = {0003-486X,1939-8980},
   MRCLASS = {53C44 (35K55 35K99 53A10 53C42 58E30 58K60)},
  MRNUMBER = {2993752},
MRREVIEWER = {Paul\ Bryan},
       DOI = {10.4007/annals.2012.175.2.7},
       URL = {https://doi.org/10.4007/annals.2012.175.2.7},
}

@article{KO82,
    author = {Kawasaki, Kyozi and Ohta, Takao},
    title = {Kinetic Drumhead Model of Interface. I},
    journal = {Progress of Theoretical Physics},
    volume = {67},
    number = {1},
    pages = {147-163},
    year = {1982},
    month = {01},
    issn = {0033-068X},
    doi = {10.1143/PTP.67.147},
    url = {https://doi.org/10.1143/PTP.67.147},
    eprint = {https://academic.oup.com/ptp/article-pdf/67/1/147/5438937/67-1-147.pdf},
}

@article {Gra87,
    AUTHOR = {Grayson, Matthew A.},
     TITLE = {The heat equation shrinks embedded plane curves to round
              points},
   JOURNAL = {J. Differential Geom.},
  FJOURNAL = {Journal of Differential Geometry},
    VOLUME = {26},
      YEAR = {1987},
    NUMBER = {2},
     PAGES = {285--314},
      ISSN = {0022-040X,1945-743X},
   MRCLASS = {53A04 (35K99 58G11)},
  MRNUMBER = {906392},
MRREVIEWER = {R.\ Osserman},
       URL = {http://projecteuclid.org/euclid.jdg/1214441371},
}

@ARTICLE{Einstein1905,
       author = {{Einstein}, A.},
        title = "{{\"U}ber die von der molekularkinetischen Theorie der W{\"a}rme geforderte Bewegung von in ruhenden Fl{\"u}ssigkeiten suspendierten Teilchen}",
      journal = {Annalen der Physik},
         year = 1905,
        month = jan,
       volume = {322},
       number = {8},
        pages = {549-560},
          doi = {10.1002/andp.19053220806},
       adsurl = {https://ui.adsabs.harvard.edu/abs/1905AnP...322..549E}
}

@article {AV22a,
    AUTHOR = {Agresti, Antonio and Veraar, Mark},
     TITLE = {Nonlinear parabolic stochastic evolution equations in critical
              spaces part {I}. {S}tochastic maximal regularity and local
              existence},
   JOURNAL = {Nonlinearity},
  FJOURNAL = {Nonlinearity},
    VOLUME = {35},
      YEAR = {2022},
    NUMBER = {8},
     PAGES = {4100--4210},
      ISSN = {0951-7715,1361-6544},
   MRCLASS = {60H15 (35K59 35R60 42B37 47D06)},
  MRNUMBER = {4459102},
       DOI = {10.1088/1361-6544/abd613},
       URL = {https://doi.org/10.1088/1361-6544/abd613},
}

@article {AV22b,
    AUTHOR = {Agresti, Antonio and Veraar, Mark},
     TITLE = {Nonlinear parabolic stochastic evolution equations in critical
              spaces part {II}: {B}low-up criteria and instataneous
              regularization},
   JOURNAL = {J. Evol. Equ.},
  FJOURNAL = {Journal of Evolution Equations},
    VOLUME = {22},
      YEAR = {2022},
    NUMBER = {2},
     PAGES = {Paper No. 56, 96},
      ISSN = {1424-3199,1424-3202},
   MRCLASS = {35R60 (35K59 60H15)},
  MRNUMBER = {4437443},
MRREVIEWER = {Le\ Chen},
       DOI = {10.1007/s00028-022-00786-7},
       URL = {https://doi.org/10.1007/s00028-022-00786-7},
}

\end{document}